\documentclass[matprg,final]{svjour}
\smartqed  
\usepackage{biograph}        
\usepackage{amssymb,amsmath}
\usepackage{graphicx}
\usepackage{times}
\usepackage{hyperref}

\def\RR{\mathbb{R}}
\def\SS{\mathbb{S}}

\def\NN{\mathbb{N}}

\def\Tr{\mbox{\rm Tr\,}}
\def\Vec{\mbox{\rm vec\,}}

\newtheorem{Algorithm}{Algorithm}

\title{PENLAB: A MATLAB solver for nonlinear semidefinite optimization%
\thanks{The research of MK was partly supported
by the Grant Agency of the Czech Republic through project GAP201-12-0671.}}
\author{Jan Fiala\ \and Michal Ko\v{c}vara\ \and Michael Stingl}
\institute{
Jan Fiala\at The Numerical Algorithms Group Ltd, Wilkinson House,
Jordan Hill Road, Oxford, OX2 8DR, UK, \email{jan@nag.co.uk} \and
Michal Ko\v{c}vara\at School of Mathematics, University of
Birmingham, Birmingham B15 2TT, UK and Institute of Information Theory
and Automation, Academy of Sciences of the Czech Republic, Pod
vod\'arenskou v\v{e}\v{z}\'{\i}~4, 18208 Praha 8,
Czech Republic, \email{m.kocvara@bham.ac.uk}
   \and Michael Stingl\at
Applied Mathematics II, University of
Erlangen-Nuremberg, N\"a{}gelsbachstr. 49b, 91052 Erlangen, Germany,
\email{stingl@am.uni-erlangen.de}}
\date{Received: date / Revised version: date}

\begin{document}
\maketitle

\begin{abstract}
PENLAB is an open source software package for nonlinear optimization,
linear and nonlinear semidefinite optimization and any combination of
these. It is written entirely in MATLAB. PENLAB is a young brother of
our code PENNON \cite{pennon} and of a new implementation from NAG
\cite{naglib}: it can solve the same classes of problems and uses the same
algorithm. Unlike PENNON, PENLAB is open source and allows the user not
only to solve problems but to modify various parts of the algorithm. As
such, PENLAB is particularly suitable for teaching and research
purposes and for testing new algorithmic ideas.

In this article, after a brief presentation of the underlying
algorithm, we focus on practical use of the solver, both for general
problem classes and for specific practical problems.
\end{abstract}

\section{Introduction}
Many problems in various scientific disciplines, as well as many
industrial problems lead to (or can be advantageously formulated) as
nonlinear optimization problems with semidefinite constraints. These
problems were, until recently, considered numerically unsolvable, and
researchers were looking for other formulations of their problem that
often lead only to approximation (good or bad) of the true solution.
This was our main motivation for the development of PENNON
\cite{pennon}, a code for nonlinear optimization problems with matrix
variables and matrix inequality constraints.

Apart from PENNON, other concepts for the solution of nonlinear
semidefinite programs are suggested in literature; see
\cite{sun-sun-zhang} for a discussion on the classic augmented
Lagrangian method applied to nonlinear semidefinite programs,
\cite{correa2004global,fares,freund-jarre-vogelbusch} for sequential
semidefinite programming algorithms and \cite{kanzow-nagel-newt} for a
smoothing type algorithm. However, to our best knowledge, none of these
algorithmic concepts lead to a publicly available code yet.

In this article, we present PENLAB, a younger brother of PENNON and
a new implementation from NAG. PENLAB can solve the same classes of
problems, uses the same algorithm and its behaviour is very similar.
However, its performance is relatively limited in comparison to \cite{pennon}
and \cite{naglib}, due to MATLAB implementation. On the other hand,
PENLAB is open source and allows the user not only to solve problems
but to modify various parts of the algorithm. As such, PENLAB is
particularly suitable for teaching and research purposes and for
testing new algorithmic ideas.

After a brief presentation of the underlying algorithm, we focus on
practical use of the solver, both for general problem classes and for
specific practical problems, namely, the nearest correlation matrix
problem with constraints on condition number, the truss topology
problem with global stability constraint and the static output
feedback problem. More applications of nonlinear semidefinite
programming problems can be found, for instance, in
\cite{annad,kanno,leibfritz-volkwein}.

PENLAB is distributed under GNU GPL license and can be downloaded from
{\tt http://web.mat.bham.ac.uk/kocvara/penlab}.

We use standard notation: Matrices are denoted by capital letters
($A,B,X,\ldots$) and their elements by the corresponding small-case
letters ($a_{ij}, b_{ij}, x_{ij},\ldots$). For vectors $x,y\in\RR^n$,
$\langle x,y\rangle:=\sum_{i=1}^n x_iy_i$ denotes the inner product.
$\SS^{m}$ is the space of real symmetric matrices of dimension $m\times
m$. The inner product on $\SS^{m}$ is defined by $\langle A,
B\rangle_{\SS^{m}} := \Tr (AB)$. When the dimensions of $A$ and $B$ are
known, we will often use notation $\langle A, B\rangle$, same as for
the vector inner product. Notation $A\preccurlyeq B$ for
$A,B\in\SS^{m}$ means that the matrix $B-A$ is positive semidefinite.
 If $A$ is an $m\times n$ matrix and $a_j$ its $j$-th
column, then $\Vec A$ is the $mn\times 1$ vector
$$
  \Vec A = \begin{pmatrix} a_1^T\ \ a_2^T\ \ \cdots\ \ a_n^T\end{pmatrix}^T\,.
$$
Finally, for $\Phi:\SS^m\to\SS^m$ and $X,Y\in \SS^m$, $D\Phi(X;Y)$
denotes the directional derivative of $\Phi$ with respect to $X$ in
direction $Y$.

\vfill

\section{The problem}
We intend to solve optimization problems with a nonlinear objective
subject to nonlinear inequality and equality constraints and nonlinear
matrix inequalities (NLP-SDP):
\begin{align}
& \min_{x\in \RR^n, Y_1\in\SS^{p_1},\ldots,Y_k\in\SS^{p_k}}  f(x,Y)\label{eq:nlpsdp}\\
& \begin{aligned}
  \mbox{subject to}\quad
   &g_i(x,Y)  \leq 0, \qquad &&i=1,\ldots,m_g\\
   &h_i(x,Y)  = 0, \qquad &&i=1,\ldots,m_h \\
   &{ {\cal A}_i(x,Y)\preceq 0,} \qquad &&{ i=1,\ldots,m_A}\\
   &\underline{\lambda}_i I \preceq Y_i\preceq \overline{\lambda}_i I, \qquad&&i=1,\ldots,k\,.
\end{aligned}\nonumber
\end{align}
Here
\begin{itemize}
\item $x\in\RR^n$ is the vector variable;
\item $Y_1\in\SS^{p_1},\ldots,Y_k\in\SS^{p_k}$ are the matrix
    variables, $k$ symmetric matrices of dimensions $p_1\times
    p_1,\ldots,p_k\times p_k$;
\item we denote $Y=(Y_1,\ldots,Y_k)$;
\item $f$, $g_i$ and $h_i$ are $C^2$ functions from $\RR^n\times
    \SS^{p_1}\times\ldots\times\SS^{p_k}$ to $\RR$;
\item $\underline{\lambda}_i$ and $\overline{\lambda}_i$ are the
    lower and upper bounds, respectively, on the eigenvalues of~
    $Y_i$, $i=1,\ldots,k$;
\item ${\cal A}_i(x,Y)$ are twice continuously differentiable
    nonlinear matrix operators from $\RR^n\times
    \SS^{p_1}\times\ldots\times\SS^{p_k}$ to $\SS^{p_{A_i}}$ where
    ${p_{A_i}}$, $i=1,\ldots,m_A$, are positive integers.
\end{itemize}


\section{The algorithm}
The basic algorithm used in this article is based on the nonlinear
rescaling method of Roman~Polyak \cite{polyak} and was described in
detail in \cite{pennon} and \cite{stingl}. Here we briefly recall it
and stress points that will be needed in the rest of the paper.

The algorithm is based on a choice of penalty/barrier functions
$\varphi:\RR\to\RR$ that penalize the inequality constraints and
$\Phi:\SS^p\to\SS^p$ penalizing the matrix inequalities. These
functions satisfy a number of properties (see \cite{pennon,stingl})
that guarantee that for any $\pi>0$ and $\Pi > 0$, we have
$$
  z(x) \le 0 \ \Longleftrightarrow \ \pi\varphi(z(x)/\pi) \le
0,
   \quad z\in C^2(\RR^n\to\RR)
$$
and
$$
  Z \preceq 0 \ \Longleftrightarrow \ \Pi \Phi(Z/\Pi) \preceq
0, \quad  Z\in\SS^p \,.
$$
This  means that, for any $\pi>0$, $\Pi>0$, problem (\ref{eq:nlpsdp})
has the same solution as the following ``augmented" problem
\begin{align}
& \min_{x\in\RR^n, Y_1\in\SS^{p_1},\ldots,Y_k\in\SS^{p_k}}   f(x,Y) \label{eq:nlpsdp_phi}\\
& \begin{aligned}
  \mbox{subject to}\quad
   &\varphi_{\pi}(g_i(x,Y))  \leq 0, \qquad &&i=1,\ldots,m_g\nonumber\\
   &\Phi_{\Pi}({\cal A}_i(x,Y))\preceq 0,\qquad&&i=1,\ldots,m_A\nonumber\\
   &\Phi_{\Pi}(\underline{\lambda}_i I -  Y_i)\preceq 0, \qquad&&i=1,\ldots,k\nonumber\\
   & \Phi_{\Pi}(Y_i-\overline{\lambda}_i I)\preceq 0, \qquad&&i=1,\ldots,k\nonumber\\
   &h_i(x,Y)  = 0, \qquad &&i=1,\ldots,m_h\,, \nonumber\\
\end{aligned}
\end{align}
where we have used the abbreviations $\varphi_{\pi} = \pi \varphi
(\cdot / \pi)$ and $\Phi_{\Pi} = \Pi \Phi (\cdot / \Pi)$.


\medskip
The Lagrangian of (\ref{eq:nlpsdp_phi}) can be viewed as
a (generalized) augmented Lagrangian of (\ref{eq:nlpsdp}):
\begin{multline}\label{eq:lagr}
  F(x,Y,u,\Xi,\underline{U},\overline{U},v,\pi,\Pi)\\
   =f(x,Y) + \sum_{i=1}^{m_g} u_i
\varphi_{\pi}(g_i(x,Y))
+ \sum_{i=1}^{m_A}\langle \Xi_i,\Phi_\Pi({\cal A}_i(x,Y))\rangle\\
 +
\sum_{i=1}^{k}\langle\underline{U}_i, \Phi_\Pi(\underline{\lambda}_i I
- Y_i)\rangle +
\sum_{i=1}^{k}\langle\overline{U}_i,\Phi_\Pi(Y_i-\overline{\lambda}_i
I)\rangle +v^\top h(x, Y)
  \,;
\end{multline}
here $u\in\RR^{m_g}$, $\Xi = (\Xi_1,\ldots,\Xi_{m_A}),
\Xi_i\in\SS^{p_{A_i}}$, and
$\underline{U}=(\underline{U}_1,\ldots,\underline{U}_k),
\overline{U}=(\overline{U}_1, \ldots, \overline{U}_k)$,
$\underline{U}_i,\overline{U}_i\in\SS^{p_i}$, are Lagrange multipliers
associated with the standard and the matrix inequality constraints,
respectively, and $v\in \RR^{m_h}$ is the vector of Lagrangian
multipliers associated with the equality constraints.

The algorithm combines ideas of the (exterior) penalty and (interior)
barrier methods with the augmented Lagrangian method.
\begin{Algorithm}\label{algo:1}
Let $x^1, Y^1$ and $u^1, \Xi^1, \underline{U}^1, \overline{U}^1, v^1$
be given. Let $\pi^1>0$, $\Pi^1>0$ and $\alpha^1>0$. For
$\ell=1,2,\ldots$ repeat till a stopping criterium is reached:
\begin{align*}
(i)\qquad &\mbox{Find $x^{\ell+1}$, $Y^{\ell+1}$ and $v^{\ell+1}$ such that}\\
\qquad &\qquad\|\nabla_{x,Y} F(x^{\ell+1},Y^{\ell+1},u^{\ell},\Xi^{\ell},\underline{U}^{\ell},
\overline{U}^{\ell},v^{\ell+1},\pi^{\ell},\Pi^{\ell})\| \leq \alpha^{\ell}\\
\qquad &\qquad\| h(x^{\ell+1},Y^{\ell+1}) \| \leq \alpha^{\ell}\\
(ii)\qquad &u_i^{\ell+1}  =  u_i^{\ell}\varphi_{\pi^{\ell}}'(g_i(x^{\ell+1},Y^{\ell+1})),\quad
i=1,\,\ldots,m_g\\
&\Xi_i^{\ell+1} = D_{\cal A} \Phi_{\Pi^{\ell}}({\cal A}_i(x^{\ell+1},Y^{\ell+1});
\Xi_i^{\ell}),\quad
i=1,\,\ldots,m_A\\
&\underline{U}_i^{\ell+1} = D_{\cal A} \Phi_{\Pi^{\ell}}((\underline{\lambda}_i I -  Y_i^{\ell+1});
\underline{U}_i^{\ell}),\quad
i=1,\,\ldots,k\\
&\overline{U}_i^{\ell+1} = D_{\cal A} \Phi_{\Pi^{\ell}}(( Y_i^{\ell+1}- \overline{\lambda}_i I );
\overline{U}_i^{\ell}),\quad
i=1,\,\ldots,k\\
(iii)\qquad &\pi^{\ell+1}  <  \pi^{\ell},\quad
\Pi^{\ell+1} < \Pi^{\ell},\quad \alpha^{\ell+1} < \alpha^{\ell}\,.
\end{align*}
\end{Algorithm}

In Step~(i) we attempt to find an approximate solution of the following
system (in $x, Y$ and $v$):
\begin{equation}\label{eq:KKT-eq}
\begin{aligned}
  \qquad \nabla_{x,Y} {F} (x,Y,u,\Xi,\underline{U},\overline{U},v,\pi,\Pi) &= 0 \\
 h(x,Y) &= 0\,,
\end{aligned}
\end{equation}
where the penalty parameters $\pi, \Pi$, as well as the multipliers
$u,\Xi,\underline{U},\overline{U}$ are fixed. In order to solve it, we
apply the damped Newton method. Descent directions are calculated
utilizing the MATLAB command {\tt ldl} that is based on the
factorization routine MA57, in combination with an inertia correction
strategy described in \cite{stingl}. In the forthcoming release of
PENLAB, we will also apply iterative methods, as described in
\cite{pen-iter}. The step length is derived using an augmented
Lagrangian merit function defined as
$$
  {F} (x,Y,u,\Xi,\underline{U},\overline{U},v,\pi,\Pi) + \frac{1}{2\mu}\|h(x,Y)\|_2^2
$$
along with an Armijo rule.

If there are no equality constraints in the problems, the unconstrained
minimization in Step~(i) is performed by the modified Newton method
with line-search (for details, see \cite{pennon}).

The multipliers calculated in Step~(ii) are restricted in order to
satisfy:
$$
\mu < \frac{u_i^{\ell+1}}{u_i^{\ell}} < \frac1{\mu}
$$
with some positive ${\mu} \leq 1$; by default, ${\mu} = 0.3$.
A similar restriction procedure can be applied to the matrix
multipliers $\underline{U}^{\ell+1}, \overline{U}^{\ell+1}$ and $\Xi$;
see again \cite{pennon} for details.

The penalty parameters $\pi, \Pi$ in Step~(iii) are updated by some
constant factor dependent on the initial penalty parameters $\pi^1,
\Pi^1$. The update process is stopped when $\pi_{eps}$ (by default
$10^{-6}$) is reached.

Algorithm~\ref{algo:1} is stopped when a criterion based on the
KKT error is satisfied and both of the inequalities holds:
\begin{eqnarray*}
\frac{|f(x^{\ell},Y^{\ell}) - F(x^{\ell},Y^{\ell},u^{\ell},\Xi^{\ell},\underline{U}^{\ell},\overline{U}^{\ell},v^{\ell},\pi^{\ell},\Pi^{\ell})|}
{1+|f(x^{\ell},Y^{\ell})|} &<& \epsilon\\
\frac{|f(x^{\ell},Y^{\ell}) - f(x^{\ell-1},Y^{\ell-1})|}{1+|f(x^{\ell},Y^{\ell})|} &<& \epsilon\,,
\end{eqnarray*}
where $\epsilon$ is by default $10^{-6}$.

\subsection{{Choice of $\varphi$ and $\Phi$}}\label{sec:hess}
To treat the standard NLP constraints, we use the penalty/barrier
function proposed by Ben-Tal and Zibulevsky \cite{ben-tal-zibulevsky}:
\begin{equation}
\varphi_{\bar{\tau}} (\tau) = \left \{
   \begin{aligned}
      &\tau + \frac{1}{2} \, \tau^2 &\mbox{if~}& \tau \geq \bar{\tau} \\
      &- (1+ \bar{\tau})^2 \log \left ( \frac{1+ 2 \bar{\tau} -\tau}
	 {1 + \bar{\tau}} \right)
      + \bar{\tau} + \frac{1}{2} \bar{\tau}^2 \  &\mbox{if~}& \tau < \bar{\tau}  \,;
    \end{aligned} \right .
\label{eq:phi}
\end{equation}
by default, $\bar{\tau} = - \frac{1}{2}$.

The penalty function $\Phi_\Pi$ of our choice is defined as follows
(here, for simplicity, we omit the variable $Y$):
\begin{equation}\label{eq:pen}
\Phi_\Pi({\cal A}(x))  = -\Pi^2({\cal A}(x) - \Pi I)^{-1} - \Pi I \,.
\end{equation}
The advantage of this choice is that it gives closed formulas for the
first and second derivatives of $\Phi_\Pi$. Defining
\begin{equation}\label{eq:Z}
  {\cal Z}(x) = -({\cal A}(x) - \Pi I)^{-1}
\end{equation}
we have (see \cite{pennon}):
\begin{align*}
\frac{\partial}{\partial x_i} \Phi_\Pi({\cal A}(x))& =
\Pi^2{\cal Z}(x) \frac{\partial{\cal A}(x)}{\partial x_i} {\cal Z}(x)
\label{eq:der1} \\
\frac{\partial^2}{\partial x_i\partial x_j} \Phi_\Pi({\cal
A}(x)) & = \Pi^2{\cal Z}(x) \left(\frac{\partial{\cal A}(x)}{\partial
x_i}
{\cal Z}(x) \frac{\partial{\cal A}(x)}{\partial x_j} +
 \frac{\partial^2{\cal A}(x)}{\partial x_i\partial x_j}
  \right.\nonumber\\
&\left.\phantom{\Pi^2{\cal Z}(x)}\qquad\ \
+ \frac{\partial{\cal A}(x)}{\partial x_j}
{\cal Z}(x) \frac{\partial{\cal A}(x)}{\partial x_i}\right){\cal
Z}(x)\,.&\nonumber
\end{align*}

\subsection{Strictly feasible constraints} In certain applications,
some of the bound constraints must remain strictly feasible for all
iterations because, for instance, the objective function may be
undefined at infeasible points (see examples in
Section~\ref{ex:truss}). To be able to solve such problems, we treat
these inequalities by a classic barrier function. In case of matrix
variable inequalities, we
split $Y$ in non-strictly feasible matrix variables $Y_1$ and strictly
feasible matrix variables $Y_2$, respectively, and define the augmented
Lagrangian
\begin{equation}\label{eq:lagr2}
\widetilde{F}(x,Y_1,Y_2,u,\Xi,\underline{U},\overline{U},v,\pi,\Pi,\kappa) =
F(x,Y_1,u,\Xi,\underline{U},\overline{U},v,\pi,\Pi) + \kappa \Phi_{\rm bar}(Y_2),
\end{equation}
where $\Phi_{\rm bar}$ can be defined, for example for the constraint
$Y_2\succeq 0$, by
$$\Phi_{\rm bar}(Y_2) = -\log\det(Y_2).$$
Strictly feasible variables $x$ are treated in a similar manner.
Note that, while the penalty parameter $\pi$ may be constant from a
certain index $\bar{\ell}$ (see again \cite{stingl} for details), the
barrier parameter $\kappa$ is required to tend to zero with increasing
$\ell$.


\section{The code}
PENLAB is a free open-source MATLAB implementation of the algorithm
described above. The main attention was given to clarity of the code
rather than tweaks to improve its performance. This should allow users
to better understand the code and encourage them to edit and develop
the algorithm further. The code is written entirely in MATLAB with an
exception of two mex-functions that handles the computationally most
intense task of evaluating the second derivative of the Augmented
Lagrangian and a sum of multiple sparse matrices (a slower non-mex
alternative is provided as well). The
solver is implemented as a MATLAB handle class and thus it should be
supported on all MATLAB versions starting from R2008a.

PENLAB is distributed under GNU GPL license and can be downloaded from
{\tt http://web.mat.bham.ac.uk/kocvara/penlab}. The distribution
package includes the full source code and precompiled mex-functions,
PENLAB User's Guide and also an internal (programmer's) documentation
which can be generated from the source code. Many examples provided in
the package show various ways of calling PENLAB and handling NLP-SDP
problems.

\subsection{Usage}
The source code is divided between a class \verb|penlab| which
implements Algorithm 1 and handles generic NLP-SDP problems similar to
formulation (\ref{eq:nlpsdp}) and interface routines providing various
specialized inputs to the solver.
Some of these are described in Section~\ref{sec:modules}.

The user needs to prepare a MATLAB structure (here called \verb|penm|)
which describes the problem parameters, such as number of variables,
number of constraints, lower and upper bounds, etc. Some of the fields
are shown in Table~\ref{tab:7}, for a complete list see the PENLAB
User's Guide. The structure is passed to \verb|penlab| which returns
the initialized problem instance:
\begin{verbatim}
>> problem = penlab(penm);
\end{verbatim}
The solver might be invoked and results retrieved, for example, by
calling
\begin{verbatim}
>> problem.solve()
>> problem.x
\end{verbatim}

The point \verb|x| or option settings might be changed and
the solver invoked again. The whole object can be cleared from the
memory using
\begin{verbatim}
>> clear problem;
\end{verbatim}

\begin{table}[htbp]
  \caption{Selection of fields of the MATLAB structure {\tt penm} used to
initialize PENLAB object. Full list is available in PENLAB
User's Guide.}
  \begin{tabular*}{\hsize}{@{\extracolsep{\fill}}ll}
\hline
field name  & meaning \\
\hline
  Nx & dimension of vector $x$ \\
  NY & number of matrix variables $Y$\\
  Y  & cell array of length NY with a nonzero pattern of each of the matrix variables\\
  lbY & NY lower bounds on matrix variables (in spectral sense) \\
  ubY & NY upper bounds on matrix variables (in spectral sense)  \\
  NANLN &  number of nonlinear matrix constraints     \\
  NALIN &  number of linear matrix constraints       \\
  lbA & lower bounds on all matrix constraints\\
  ubA & upper bounds on all matrix constraints\\
\hline
\end{tabular*}   \label{tab:7}
\end{table}

\subsection{Callback functions} The principal philosophy of the code is
similar to many other optimization codes---we use callback functions
(provided by the user) to compute function values and derivatives of
all involved functions.

For a generic problem, the user must define nine MATLAB callback
functions: {\tt objfun}, {\tt objgrad}, {\tt objhess}, {\tt confun},
{\tt congrad}, {\tt conhess}, {\tt mconfun}, {\tt mcongrad}, {\tt
mconhess} for function value, gradient, and Hessian of the objective
function, (standard) constraints and matrix constraint. If one
constraint type is not present, the corresponding callbacks need not be
defined. Let us just show the parameters of the most complex callbacks
for the matrix constraints:
\begin{verbatim}
function [Ak, userdata] = mconfun(x,Y,k,userdata)
function [dAki,userdata] = mcongrad(x,Y,k,i,userdata)
function [ddAkij, userdata] = mconhess(x,Y,k,i,j,userdata)
\end{verbatim}
Here $x,Y$ are the current values of the (vector and matrix) variables.
Parameter $k$ stands for the constraint number. Because every element
of the gradient and the Hessian of a matrix function is a matrix, we
compute them (the gradient and the Hessian) element-wise (parameters
$i,j$).  The outputs {\tt
Ak,dAki,ddAkij} are symmetric matrices saved in sparse MATLAB format.

Finally, {\tt userdata} is a MATLAB structure passed through
all callbacks for user's convenience and may contain any
additional data needed for the evaluations.
It is unchanged by the algorithm itself but it can be modified in the
callbacks by user.
For instance, some time-consuming computation that depends
on $x,Y,k$ but is independent of $i$ can be performed only for $i=1$,
the result stored in {\tt userdata} and recalled for any $i>1$ (see,
e.g., Section~\ref{ex:truss}, example Truss Design with Buckling Constraint).

\subsection{Mex files}

Despite our intentions to use only pure Matlab code, two routines were
identified to cause a significant slow-down and therefore their m-files
were substituted with equivalent mex-files. The first one computes
linear combination of a set of sparse matrices, e.g., when evaluating
${\cal A}_i(x)$ for polynomial matrix inequalities, and is based on
ideas from \cite{davis}. The second one evaluates matrix inequality
contributions to the Hessian of the augmented Lagrangian
(\ref{eq:lagr}) when using penalty function (\ref{eq:pen}).

The latter case reduces to computing $z_{\ell} = \langle TA_kU,\,
A_{\ell}\rangle$ for $\ell=k,\ldots,n$ where $T, U \in \SS^m$ are dense
and $A_{\ell} \in \SS^m$ are sparse with potentially highly varying
densities. Such expressions soon become challenging for nontrivial $m$
and can easily dominate the whole Algorithm~\ref{algo:1}. Note that the
problem is common even in primal-dual interior point methods for SDPs
and have been studied in \cite{fujisawa-kojima-nakata}. We developed a
relatively simple strategy which can be viewed as an evolution of the
three computational formulae presented in \cite{fujisawa-kojima-nakata}
and offers a minimal number of multiplications while keeping very
modest memory requirements. We refer to it as a \emph{look-ahead
strategy with caching}. It can be described as follows:

\begin{Algorithm}\label{algo:trace}
Precompute a set ${\cal J}$ of all nonempty columns across all
$A_{\ell}, {\ell}=k,\ldots,n$ and a set ${\cal I}$ of nonempty rows of
$A_k$ \emph{(look-ahead)}. Reset flag vector $c\leftarrow 0$, set $z=0$
and $v=w=0$. For each $j \in {\cal J}$ perform:
\begin{enumerate}
\item compute selected elements of the $j$-th column of $A_kU$, i.e.,\\
    $v_i = \sum_{\alpha=1}^m(A_k)_{i \alpha} U_{\alpha j}$ for $i \in {\cal I}$,
\item for each $A_{\ell}$ with nonempty $j$-th column go through
    its nonzero elements $(A_{\ell})_{ij}$ and
  \begin{enumerate}
    \item if $c_i<j$ compute $w_i = \sum_{\alpha \in {\cal I}}
        T_{i\alpha}v_{\alpha}$ and set $c_i \leftarrow j$
        \emph{(caching)},
    \item update trace, i.e., $z_{\ell} = z_{\ell} + w_i(A_{\ell})_{ij}$.
  \end{enumerate}
\end{enumerate}
\end{Algorithm}


\section{Gradients and Hessians of matrix valued functions}

There are several concepts of derivatives of matrix functions; they,
however, only differ in the ordering of the elements of the resulting
``differential''. In PENLAB, we use the following definitions of the
gradient and Hessian of matrix valued functions.
\begin{definition}\label{def:1}
Let $F$ be a differentiable $m\times n$ real matrix function of an
$p\times q$ matrix of real variables $X$. The $(i,j)$-th element of the
\emph{gradient} of $F$ at $X$ is the $m\times n$ matrix
\begin{equation}\label{eq:a041}
\left[\nabla F(X)\right]_{ij} :=
\frac{\partial F(X)}{\partial x_{ij}},
  \qquad i=1,\ldots,p,\ j=1,\ldots,q
\,.
\end{equation}
\end{definition}
\begin{definition}\label{def:2}
Let $F$ be a twice differentiable $m\times n$ real matrix function of
an $p\times q$ matrix of real variables $X$. The $(ij,k\ell)$-th
element of the \emph{Hessian} of $F$ at $X$ is the $m\times n$ matrix
\begin{equation}\label{eq:a040}
\left[\nabla^2 F(X)\right]_{ij,k\ell} :=
\frac{\partial^2 F(X)}{\partial x_{ij}\partial x_{kl}} ,
  \qquad i,k=1,\ldots,p,\ j,\ell=1,\ldots,q
\,.
\end{equation}
\end{definition}
In other words, for every pair of variables $x_{ij},\ x_{k\ell}$,
elements of $X$, the second partial derivative of $F(X)$ with respect
to these variables is the $m\times n$ matrix $\frac{\partial^2
F(X)}{\partial x_{ij}\partial x_{k\ell}}$.

How to compute these derivatives, i.e., how to define the callback
functions? In Appendix A, we summarize basic formulas for the
computation of derivatives of scalar and matrix valued functions of
matrices.

For low-dimensional problems, the user can utilize MATLAB's Symbolic
Toolbox. For instance, for $F(X)=XX$, the commands
\begin{verbatim}
>> A=sym('X',[2,2]);
>> J=jacobian(X*X,X(:));
>> H=jacobian(J,X(:));
\end{verbatim}
generate arrays $J$ and $H$ such that the $i$-th column of $J$ is the
vectorized $i$-th element of the gradient of $F(X)$; similarly, the
$k$-th column of $H$, $k=(i-1)n^2+j$ for $i,j=1,\ldots,n^2$ is the
vectorized $(i,j)$-th element of the Hessian of $F(X)$. Clearly, the
dimension of the matrix variable is fixed and for a different dimension
we have to generate new formulas. Unfortunately, this approach is
useless for higher dimensional matrices (the user is invited to use the
above commands for $F(X)=X^{-1}$ with $X\in\SS^5$ to see the
difficulties). However, one can always use symbolic computation to
check validity of general dimension independent formulas on small
dimensional problems.

\section{Pre-programmed interfaces}\label{sec:modules}
PENLAB distribution contains several pre-programmed interfaces for
standard optimization problems with standard inputs. For these
problems, the user does not have to create the \verb|penm| object, nor
the callback functions.

\subsection{Nonlinear optimization with AMPL input}
PENLAB can read optimization problems that are defined in and processed
by AMPL \cite{ampl}. AMPL contains routines for automatic
differentiation, hence the gradients and Hessians in the callbacks
reduce to calls to appropriate AMPL routines.

Assume that nonlinear optimization problem is processed by AMPL, so
that we have the corresponding \verb|.nl| file, for instance
\verb|chain.nl|, stored in directory \verb|datafiles|. All the user has
to do to solve the problem is to call the following three commands:

\begin{verbatim}
>> penm = nlp_define('datafiles/chain100.nl');
>> problem = penlab(penm);
>> problem.solve();
\end{verbatim}

\subsection{Linear semidefinite programming}
Assume that the data of a linear SDP problem is stored in a MATLAB
structure \verb|sdpdata|. Alternatively, such a structure can be
created by the user from SDPA input file \cite{sdpa}. For instance, to
read problem \verb|arch0.dat-s| stored in directory \verb|datafiles|,
call
\begin{verbatim}
>> sdpdata = readsdpa('datafiles/control1.dat-s');
\end{verbatim}
To solve the problem by PENLAB, the user just has to call the following
sequence of commands:
\begin{verbatim}
>> penm = sdp_define(sdpdata);
>> problem = penlab(penm);
>> problem.solve();
\end{verbatim}

\subsection{Bilinear matrix inequalities}
We want to solve an optimization problem with  quadratic objective and
constraints in the form of bilinear matrix inequalities:
\begin{align}\label{eq:bmiproblem}
 &\min_{x\in\RR^n} \frac{1}{2} x^T H x + c^T x \\
 & \begin{aligned}
  \mbox{subject to}\quad
  &b_{\rm low}\leq Bx\leq {b_{\rm up}}\nonumber\\
 &Q^i_0 + \sum_{k=1}^{n} x_k Q^i_k + \sum_{k=1}^{n}\sum_{\ell=1}^{n} x_k x_\ell Q^i_{k\ell}
 \succcurlyeq 0, \quad i=1,\ldots,m\,. \nonumber
 \end{aligned}
\end{align}

The problem data should be stored in a simple format explained in
PENLAB User's Guide.
All the user has to do to solve the problem
is to call the following sequence of commands:

\begin{verbatim}
>> load datafiles/bmi_example;
>> penm = bmi_define(bmidata);
>> problem = penlab(penm);
>> problem.solve();
\end{verbatim}

\subsection{Polynomial matrix inequalities}
We want to solve an optimization problem with constraints in the form
of polynomial matrix inequalities:
\begin{align}\label{eq:pmiproblem}
 &\min_{x\in\RR^n} \frac{1}{2} x^T H x + c^T x \\
 & \begin{aligned}
  \mbox{subject to}\quad
  &b_{\rm low}\leq Bx\leq {b_{\rm up}}\nonumber\\
 &{\cal A}_i(x)\succcurlyeq 0,\quad i=1,\ldots,m \nonumber
 \end{aligned}
\end{align}
with
$$
 {\cal A}_i(x) = \sum_j  x^{(\kappa^i(j))} Q^i_j
$$
where $\kappa^i(j)$ is a multi-index of the $i$-th constraint with
possibly repeated entries and $x^{(\kappa^i(j))}$ is a product of
elements with indices in $\kappa^i(j)$.

For example, for $${\cal A}(x) = Q_1 + x_1 x_3 Q_2 + x_2 x_4^3 Q_3$$
the multi-indices are $\kappa(1) = \{0\}$ ($Q_1$ is an absolute term),
$\kappa(2) = \{1,3\}$ and $\kappa(3) = \{2,4,4,4\}$.

Assuming now that the problem is stored in a structure \verb|pmidata|
(as explained in PENLAB User's Guide), the user just has to call the
following sequence of commands:

\begin{verbatim}
>> load datafiles/pmi_example;
>> penm = pmi_define(pmidata);
>> problem = penlab(penm);
>> problem.solve();
\end{verbatim}

\section{Examples}
All MATLAB programs and data related to the examples in this section
can be found in directories \verb|examples| and \verb|applications|
of the PENLAB distribution.
\subsection{Correlation matrix with the constrained condition number}\label{ex:cond}
We consider the problem of finding the nearest correlation matrix
(\cite{higham}):
\begin{align}
&\min_X \sum_{i,j=1}^n (X_{ij}-H_{ij})^2\label{corr1}\\
&\mbox{subject to}\nonumber\\
&\qquad X_{ii} = 1,\quad i=1,\ldots,n\nonumber\\
&\qquad X\succeq 0\,.\nonumber
\end{align}
In addition to this standard setting of the problem, let us bound the
condition number of the nearest correlation matrix by adding the
constraint
$$
  \mbox{cond}(X) = \kappa \,.
$$
We can formulate this constraint as
\begin{align}
  I\preceq\widetilde{X}\preceq \kappa I
\end{align}
the variable transformation
$$
  \widetilde{X} = \zeta X\,.
$$
After the change of variables, and with the new constraint, the problem
of finding the nearest correlation matrix with a given condition number
reads as follows:
\begin{align}
&\min_{\zeta,\widetilde{X}} \sum_{i,j=1}^n
(\frac{1}{\zeta}\widetilde{X}_{ij}-H_{ij})^2\label{corr_cond}\\
&\mbox{subject to}\nonumber\\
&\qquad  \widetilde{X}_{ii} -\zeta = 0,\quad i=1,\ldots,n\nonumber\\
&\qquad I\preceq\widetilde{X}\preceq \kappa I\nonumber
\end{align}
The new problem now has the NLP-SDP structure of (\ref{eq:nlpsdp}).

We will consider an example based on a practical application from
finances; see \cite{werner-schoettle}. Assume that we are given a
$5\times 5$ correlation matrix. We now add a new asset class, that
means, we add one row and column to this matrix. The new data is based
on a different frequency than the original part of the matrix, which
means that the new matrix is no longer positive definite:
$$
  H_{\rm ext} = \begin{pmatrix}1 &-0.44& -0.20 &0.81& -0.46& -0.05\\
    -0.44& 1 &0.87& -0.38& 0.81 & -0.58\\
    -0.20 &.87 &1& -0.17& 0.65& -0.56\\
    0.81 &-0.38& -0.17& 1& -0.37& -0.15\\
    -0.46& 0.81& 0.65& -0.37& 1& -0.08\\
    -0.05&-0.58&-0.56&-0.15&0.08&1
    \end{pmatrix}\,.
$$

When solving problem (\ref{corr_cond}) by PENLAB with $\kappa=10$, we
get the solution after 11 outer and 37 inner iterations. The optimal
value of $\zeta$ is $3.4886$ and, after the back substitution $X =
\frac{1}{\zeta}\widetilde{X}$, we get the nearest correlation matrix
\begin{verbatim}
X =
   1.0000  -0.3775  -0.2230   0.7098  -0.4272  -0.0704
  -0.3775   1.0000   0.6930  -0.3155   0.5998  -0.4218
  -0.2230   0.6930   1.0000  -0.1546   0.5523  -0.4914
   0.7098  -0.3155  -0.1546   1.0000  -0.3857  -0.1294
  -0.4272   0.5998   0.5523  -0.3857   1.0000  -0.0576
  -0.0704  -0.4218  -0.4914  -0.1294  -0.0576   1.0000
\end{verbatim}
with eigenvalues
\begin{verbatim}
eigenvals =
   0.2866   0.2866   0.2867   0.6717   1.6019   2.8664
\end{verbatim}
and the condition number equal to 10, indeed.

\paragraph{Gradients and Hessians}
What are the first and second partial derivatives of functions involved
in problem (\ref{corr_cond})? The constraint is linear, so the answer
is trivial here, and we can only concentrate on the objective function
\begin{equation}\label{eq:f}
  f(z,\widetilde{X}):=\sum_{i,j=1}^n (z\widetilde{X}_{ij}-H_{ij})^2
  = \langle z\widetilde{X}-H, z\widetilde{X}-H \rangle\,,
\end{equation}
where, for convenience, we introduced a variable $z=\frac{1}{\zeta}$.

\begin{theorem} Let
$x_{ij}$ and $h_{ij}$, $i,j=1,\ldots,n$ be elements of $\widetilde{X}$
and $H$, respectively. For the function $f$ defined in (\ref{eq:f}) we
have the following partial derivatives:
\begin{itemize}
\item[(i)] $\nabla_{\!z}\, f(z,\widetilde{X}) = 2\langle
    \widetilde{X}, z\widetilde{X}-H\rangle$
\item[(ii)]
    $\left[\nabla_{\!\widetilde{X}}f(z,\widetilde{X})\right]_{ij} =
    2z(zx_{ij}-h_{ij})$,\quad $i,j=1,\ldots ,n$
\item[(iii)] $\nabla^2_{\!z,z}\, f(z,\widetilde{X}) = 2\langle
    \widetilde{X},\widetilde{X}\rangle$
\item[(iv)] $\left[\nabla^2_{\!z,\widetilde{X}}\,
    f(z,\widetilde{X})\right]_{ij} =
    \left[\nabla^2_{\!\widetilde{X},z}\,
    f(z,\widetilde{X})\right]_{ij} = 4zx_{ij} - 2h_{ij}$,\quad
    $i,j=1,\ldots ,n$
\item[(v)] $\left[\nabla^2_{\!\widetilde{X},\widetilde{X}}\,
    f(z,\widetilde{X})\right]_{ij,k\ell} = 2z^2$\quad for $i=k,\
    j=\ell$ and zero otherwise ($i,j,k,\ell=1,\ldots ,n$)\,.
\end{itemize}
\end{theorem}
The proof follows directly from formulas in Appendix~A.

\paragraph{PENLAB distribution}
This problem is stored in directory \verb|applications/CorrMat| of the PENLAB
distribution. To solve the above example and to see the resulting
eigenvalues of $X$, run in its directory
\begin{verbatim}
>> penm = corr_define;
>> problem = penlab(penm);
>> problem.solve();
>> eig(problem.Y{1}*problem.x)
\end{verbatim}

%
%

\subsection{Truss topology optimization with stability constraints}\label{ex:truss}
In truss optimization we want to design a pin-jointed framework
consisting of $m$ slender bars of constant mechanical properties
characterized by their Young's modulus $E$. We will consider trusses in
a $d$-dimensional space, where $d=2$ or $d=3$. The bars are jointed at
$\tilde{n}$ nodes. The system is under load, i.e., forces
$f_j\in\RR^{d}$ are acting at some nodes $j$. They are aggregated in a
vector $f$, where we put $f_j=0$ for nodes that are not under load.
This external load is transmitted along the bars causing displacements
of the nodes that make up the displacement vector $u$. Let $p$ be the
number of fixed nodal coordinates, i.e., the number of components with
prescribed discrete homogeneous Dirichlet boundary condition. We omit
these fixed components from the problem formulation reducing thus the
dimension of $u$ to
$$
n=d\,\cdot\,\tilde{n} - p .
$$
Analogously, the external load $f$ is considered as a vector in
$\RR^n$.

The design variables in the system are the bar volumes
$x_1,\ldots,x_m$. Typically, we want to minimize the weight of the
truss. We assume to have a unique material (and thus density) for all
bars, so this is equivalent to minimizing the volume of the truss,
i.e., $\sum_{i=1}^m x_i$. The optimal truss should satisfy mechanical
equilibrium conditions:
\begin{equation}
  K(x)u=f \,;
  \label{eq:3b1}
\end{equation}
here
\begin{equation}
  K(x):= \sum\limits_{i=1}^m x_iK_i,\quad K_i=\frac{E_i}{\ell_i^2}\gamma_i\gamma_i^{\top}
  \label{KO5eq:1}
\end{equation}
is the so-called stiffness matrix, $E_i$ the Young modulus of the $i$th
bar, $\ell_i$ its length and $\gamma_i$ the $n-$vector of direction
cosines.

We further introduce the compliance of the truss $f^{\top}u$ that
indirectly measures the stiffness of the structure under the force $f$
and impose the constraints
$$
  f^{\top}u \leq \gamma\,.
$$
This constraint, together with the equilibrium conditions, can be
formulated as a single linear matrix inequality (\cite{buck})
$$
\begin{pmatrix} K(x) & f\\ f^T &\gamma\end{pmatrix}
     \succeq 0\,.
$$

The minimum volume single-load truss topology optimization problem can
then be formulated as a linear semidefinite program:
\begin{align}
  &\min_{x\in\RR^m} \sum_{i=1}^m x_i\label{minvolc}\\
  &\mbox{subject to}\nonumber\\
  &\qquad \begin{pmatrix} K(x) & f\\ f^T &\gamma\end{pmatrix}
     \succeq 0 \nonumber\\
  &\qquad x_i\geq 0,\quad i=1,\ldots,m\,.\nonumber
\end{align}

We further consider the constraint on the global stability of the
truss. The meaning of the constraint is to avoid global buckling of the
optimal structure. We consider the simplest formulation of the buckling
constraint based on the so-called linear buckling assumption
\cite{buck}. As in the case of free vibrations, we need to constrain
eigenvalues of the generalized eigenvalue problem
\begin{equation}\label{eq:buckEVP}
  K(x) w = \lambda {G}(x) w \,,
\end{equation}
in particular, we require that all eigenvalues of (\ref{eq:buckEVP})
lie outside the interval [0,1]. The so-called geometry stiffness matrix
${G}(x)$ depends, this time, nonlinearly on the design variable $x$:
\begin{equation}\label{eq:G}
  {G}(x) =  \sum_{i=1}^m {G}_i(x), \qquad
  {G}_i(x) = \frac{E x_i}{\ell_i^d} (\gamma_i^{\top} K(x)^{-1}f)
  (\delta_i\delta_i^{\top}+\eta_i\eta_i^{\top}).
\end{equation}
Vectors $\delta,\eta$ are chosen so that $\gamma,\delta,\eta$ are
mutually orthogonal. (The presented formula is for $d=3$. In the
two-dimensional setting the vector $\eta$ is not present.) To simplify
the notation, we denote $$\Delta_i = \delta_i\delta^T_i +
\eta_i\eta^T_i\,.$$

It was shown in \cite{buck} that the eigenvalue constraint can be
equivalently written as a nonlinear matrix inequality
\begin{equation} \label{eq:truss_buck}
K(x)+{G}(x) \succcurlyeq 0
\end{equation}
that is now to be added to (\ref{minvolc}) to get the following
nonlinear semidefinite programming problem. Note that $x_i$ are
requested to be strictly feasible.
\begin{align}
  &\min_{x\in\RR^m} \sum_{i=1}^m x_i\label{eq:truss}\\
  &\mbox{subject to}\nonumber\\
  &\qquad \begin{pmatrix} K(x) & f\\ f^T &\gamma\end{pmatrix}
     \succeq 0
      \nonumber\\
  &\qquad K(x)+{G}(x) \succcurlyeq 0 \nonumber\\
  &\qquad x_i > 0,\quad i=1,\ldots,m\,\nonumber
\end{align}

\paragraph{Gradients and Hessians}
Let $M: \RR^m \to \RR^{n\times n}$ be a matrix valued function
assigning each vector $\xi$ a matrix $M(\xi)$. We denote by $\nabla\!_k
M$ the partial derivative of $M(\xi)$ with respect to the $k$-th
component of vector $\xi$.
\begin{lemma}[based on \cite{magnus1988matrix}]
Let $M: \RR^m \to \RR^{n\times n}$ be a symmetric matrix valued
function assigning each $\xi\in\RR^m$ a nonsingular $(n\times n)$
matrix $M(\xi)$. Then (for convenience we omit the variable $\xi$)
$$
 \nabla\!_k  M^{-1} = -M^{-1} (\nabla\!_k M)
		 M^{-1}\,.
$$
If $M$ is a linear function of $\xi$, i.e., $M(\xi) = \sum_{i=1}^m
\xi_i M_i$ with symmetric positive semidefinite $M_i, i=1,\ldots,m,$
then the above formula simplifies to
$$
 \nabla\!_k M^{-1} = -M^{-1} M_k M^{-1}\,.
$$
\end{lemma}
\begin{theorem}[\cite{buck}]
Let $G(x)$ be given as in (\ref{eq:G}). Then
$$
  [\nabla {G}\,]_k = {E\over\ell_k^3}\gamma_k^T K^{-1}f\Delta_k
    - \sum_{j=1}^m {Et_j\over\ell_j^3}\gamma_j^T K^{-1}K_k K^{-1} f
      \Delta_j
$$
and
\begin{multline*}
  \displaystyle
  [\nabla^2 {G}\,]_{\!k\ell} =
  \displaystyle
   -{E\over\ell_k^3}\gamma_k^T K^{-1}K_\ell K^{-1}f\Delta_k
   -{E\over\ell_\ell^3}\gamma_\ell^T K^{-1}K_k K^{-1}f\Delta_\ell\\
    \displaystyle
   -\sum_{j=1}^m {Et_j\over\ell_j^3}\gamma_j^T K^{-1}K_\ell K^{-1}K_k K^{-1} f
   \Delta_j\\
   -\sum_{j=1}^m {Et_j\over\ell_j^3}\gamma_j^T K^{-1}K_k K^{-1}K_\ell K^{-1} f \Delta_j.
   \end{multline*}
\end{theorem}

\begin{figure}[h]
\begin{center}
\resizebox{0.39\hsize}{!}
{\includegraphics{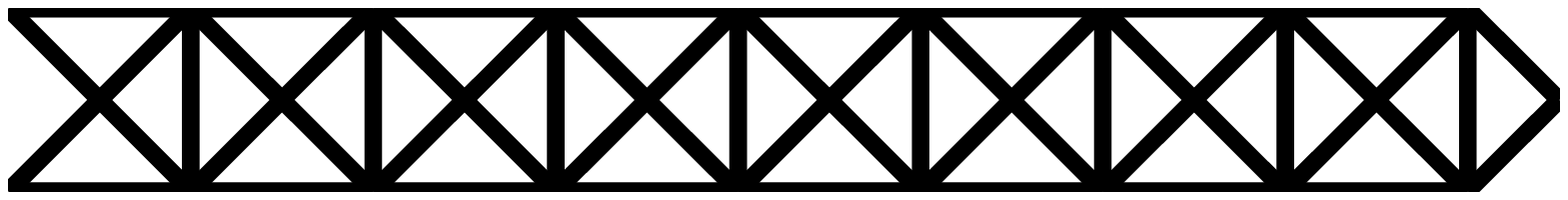}}\qquad
\resizebox{0.39\hsize}{!}
{\includegraphics{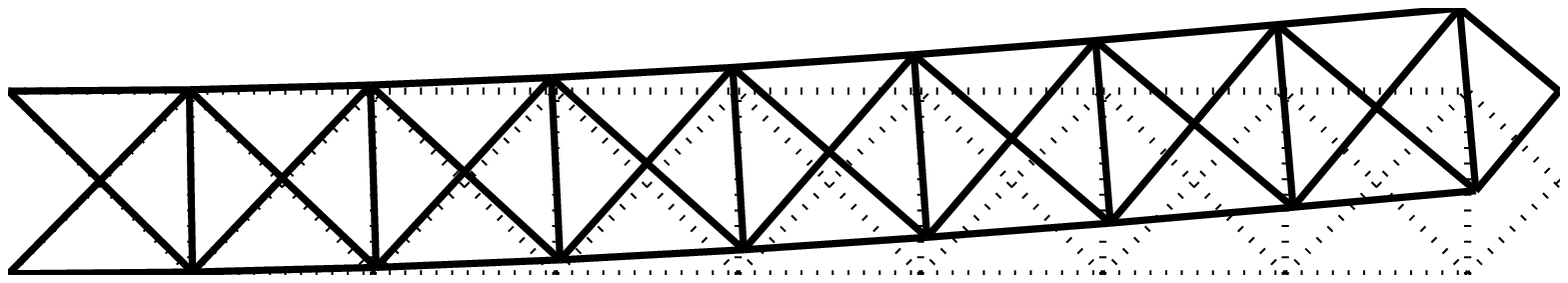}}\\[1.5em]
\resizebox{0.39\hsize}{!}
{\includegraphics{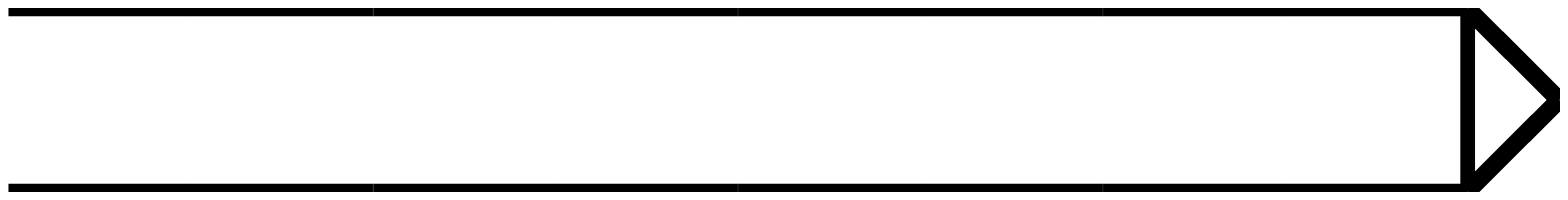}}\qquad
\resizebox{0.39\hsize}{!}
{\includegraphics{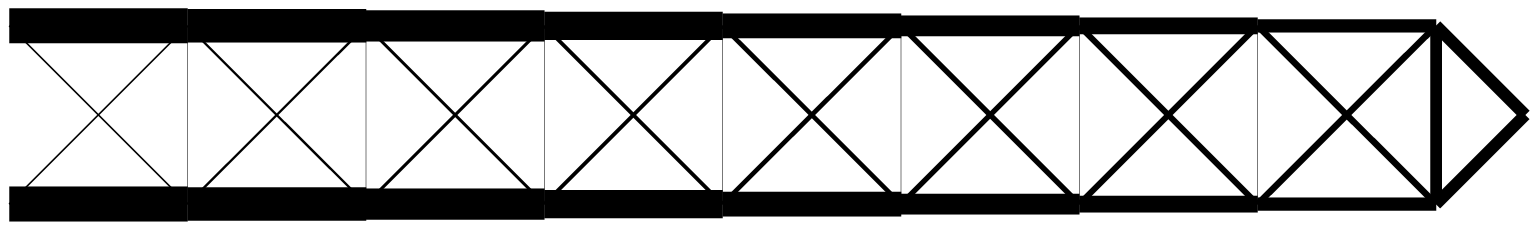}}
\end{center}
\caption{Truss optimization with stability problem:
initial truss (top-left); its buckling mode (top-right);
optimal truss without stability constraint (bottom-left);
and optimal stable truss (bottom-right)}\label{fig:tr1}
\end{figure}

\paragraph{Example} Consider the standard example of a laced column under
axial loading (example \verb|tim| in the PENLAB collection). Due to
symmetry, we only consider one half of the column, as shown in
Figure~\ref{fig:tr1}(top-peft); it has 19 nodes and 42 potential bars,
so $n=34$ and $m=42$. The column dimensions are $8.5\times 1$, the two
nodes on the left-hand side are fixed and the ``axial'' load applied at
the column tip is $(0,-10)$. The upper bound on the compliance is
chosen as $\gamma=1$.

Assume first that $x_i=0.425, i=1,\ldots,m$, i.e., the volumes of all
bars are equal and the total volume is 17.85. The values of $x_i$ were
chosen such that the truss satisfies the compliance constraint:
$f^{\top}u =0.9923\leq \gamma$. For this truss, the smallest
nonnegative eigenvalue of (\ref{eq:buckEVP}) is equal to 0.7079 and the
buckling constraint (\ref{eq:truss_buck}) is not satisfied. Figure
\ref{fig:tr1}(top-right) shows the corresponding the buckling mode
(eigenvector associated with this eigenvalue).
Let us now solve the truss optimization problem \emph{without} the
stability constraint (\ref{eq:truss}). We obtain the design shown in
Figure~\ref{fig:tr1}(bottom-left). This truss is much lighter than the
original one ($\sum\limits_{i=1}^m x_i = 9.388$), it is, however,
extremely unstable under the given load, as (\ref{eq:buckEVP}) has a
zero eigenvalue.

When solving the truss optimization problem \emph{with} the stability
constraint (\ref{eq:truss}) by PENLAB, we obtain the design shown in
Figure~\ref{fig:tr1}(bottom-right). This truss is still significantly
lighter than the original one ($\sum\limits_{i=1}^m x_i = 12.087$), but
it is now stable under the given load. To solve the nonlinear SDP
problem, PENLAB needed 18 global and 245 Newton iterations and 212
seconds of CPU time, 185 of which were spent in the Hessian evaluation
routines.

\paragraph{PENLAB distribution}
Directories \verb|applications/TTO| and \verb|applications/TTObuckling|
of the PENLAB distribution
contain the problem formulation and many examples of trusses. To solve
the above example with the buckling constraint, run
\begin{verbatim}
>> solve_ttob('GEO/tim.geo')
\end{verbatim}
in directory \verb|TTObuckling|.

\subsection{Static output feedback}

Given a linear system with $A\in\RR^{n\times n}, B\in\RR^{n\times m},
C\in\RR^{p\times n}$
\begin{align*}
\dot{x} &= Ax + Bu\\
y& = Cx
\end{align*}
we want to stabilize it by static output feedback
$
  u = Ky \,.
$ That is, we want to find a matrix $K\in\RR^{m\times p}$ such that the
eigenvalues of the closed-loop system $A+BKC$ belong to the left
half-plane.

The standard way how to treat this problem is based on the Lyapunov
stability theory. It says that $A+BKC$ has all its eigenvalues in the
open left half-plane if and only if there exists a symmetric positive
definite matrix $P$ such that
\begin{equation}\label{eq:sofbmi}
(A+BKC)^T P+P(A+BKC) \succ 0\,.
\end{equation}
Hence, by introducing the new variable, the Lyapunov matrix $P$, we can
formulate the SOF problem as a feasibility problem for the bilinear
matrix inequality (\ref{eq:sofbmi}) in variables $K$ and $P$. As
typically $n>p,m$ (often $n\gg p,m$), the Lyapunov variable dominates
here, although it is just an auxiliary variable and we do not need to
know its value at the feasible point. Hence a natural question arises
whether we can avoid the Lyapunov variable in the formulation of the
problem. The answer was given in \cite{SOF2005} and lies in the
formulation of the problem using polynomial matrix inequalities.

Let $k=\Vec K$. Define the characteristic polynomial of $A+BKC$:
$$
  q(s,k) = \det(sI-A-BKC) = \sum_{i=0}^n q_i(k)s^i\,,
$$
where $q_i(k) = \sum_\alpha q_{i\alpha}k^\alpha$ and
$\alpha\in\NN^{mp}$ are all monomial powers. The \emph{Hermite
stability criterion} says that the roots of $q(s,k)$ belong to the
stability region {\cal D} (in our case the left half-plane) if and only
if
$$
  H(q) = \sum_{i=0}^n\sum_{j=0}^n q_i(k)q_j(k) H_{ij} \succ 0 \,.
$$
Here the coefficients $H_{ij}$ depend on the stability region only
(see, e.g., \cite{ecc}). For instance, for $n=3$, we have
$$
  H(q) =  \begin{pmatrix} 2q_0q_1 & 0 & 2q_0q_3\\
                          0 & 2q_1q_2-2q_0q_3 & 0\\
                          2q_0q_3 & 0 & 2q_2q_3
          \end{pmatrix}\,.
$$
The Hermite matrix $H(q)=H(k)$ depends polynomially on $k$:
\begin{equation}\label{eq:sofpmi}
  H(k) = \sum_\alpha H_\alpha k^\alpha \succ 0
\end{equation}
where $H_\alpha = H_\alpha^T\in\RR^{n\times n}$ and $\alpha\in\NN^{mp}$
describes all monomial powers.
\begin{theorem}[\cite{SOF2005}] Matrix $K$ solves the static output feedback problem if
and only if $k=\Vec K$ satisfies the polynomial matrix inequality
(\ref{eq:sofpmi}).
\end{theorem}
In order to solve the strict feasibility problem (\ref{eq:sofpmi}), we
can solve the following optimization problem with a polynomial matrix
inequality
\begin{align}\label{eq:sofpmi1}
 &\max_{k\in\RR^{mp},\,\lambda\in\RR} \lambda-\mu\|k\|^2 \\
 &
  \mbox{subject to}\quad
  H(k)\succcurlyeq \lambda I\,. \nonumber
\end{align}
Here $\mu>0$ is a parameter that allows us to trade off between
feasibility of the PMI and a moderate norm of the matrix $K$, which is
generally desired in practice.

\paragraph{COMPlib examples}
In order to use PENLAB for the solution of SOF problems
(\ref{eq:sofpmi1}), we have developed an interface to the problem
library COMPlib \cite{complib}\footnote{The authors would like to thank
Didier Henrion, LAAS-CNRS Toulouse, for developing a substantial part
of this interface.}. Table~\ref{tab:11} presents the results of our
numerical tests. We have only solved COMPlib problems of small size,
with $n<10$ and $mp<20$. The reason for this is that our MATLAB
implementation of the interface (building the matrix $H(k)$ from
COMPlib data) is very time-consuming. For each COMPlib problem, the
table shows the degree of the matrix polynomial, problem dimensions $n$
and $mp$, the optimal $\lambda$ (the negative largest eigenvalue of the
matrix $K$), the CPU time and number of Newton iterations/linesearch
steps of PENLAB. The final column contains information about the
solution quality. ``F'' means failure of PENLAB to converge to an
optimal solution. The plus sign ``+'' means that PENLAB converged to a
solution which does not stabilize the system and "0" is used when
PENLAB converged to a solution that is on the boundary of the feasible
domain and thus not useful for stabilization.
\begin{table}[htbp]
  \caption{mmm}
  \begin{tabular*}{\hsize}{@{\extracolsep{\fill}}lcccrrrc}
\hline
Problem  & degree &  $n$ & $mp$ &  $\lambda_{\rm opt}$& CPU (sec)&iter &remark \\
\hline
 AC1 & 5 & 5 & 9 &  $-0.871\cdot 10^{0}$ & 2.2 & 27/30& \\
 AC2 & 5 & 5 & 9 &  $-0.871\cdot 10^{0}$ &  2.3 & 27/30& \\
 AC3 & 4 & 5 & 8 &  $-0.586\cdot 10^{0}$ & 1.8 & 37/48& \\
 AC4 & 2 & 4 & 2 &  $0.245\cdot 10^{-2}$ & 1.9 & 160/209& + \\
 AC6 & 4 & 7 & 8 & $-0.114\cdot 10^{4}$ & 1.2 & 22/68& \\
 AC7 & 2 & 9 & 2 &  $-0.102\cdot 10^{3}$& 0.9 & 26/91 &\\
 AC8 & 2 & 9 & 5 &  $0.116\cdot 10^{0}$ & 3.9 & 346/1276 & F \\
 AC11 & 4 & 5 & 8 &  $-0.171\cdot 10^{5}$ & 2.3 & 65/66&  \\
 AC12 & 6 & 4 & 12 &  $0.479\cdot 10^{0}$ & 12.3 & 62/73& + \\
 AC15 & 4 & 4 & 6 &  $-0.248\cdot 10^{-1}$ & 1.2 & 25/28 & \\
 AC16 & 4 & 4 & 8 &  $-0.248\cdot 10^{-1}$ & 1.2 & 23/26 & \\
 AC17 & 2 & 4 & 2 & $-0.115\cdot 10^{2}$ & 1.0 & 19/38 & \\
 HE1 & 2 & 4 & 2 & $-0.686\cdot 10^{2}$ & 1.0 & 22/22 & \\
 HE2 & 4 & 4 & 4 & $-0.268\cdot 10^{0}$ & 1.6 & 84/109 & \\
 HE5 & 4 & 8 & 8 & $0.131\cdot 10^{2}$ & 1.9 & 32/37 & + \\
 REA1 & 4 & 4 & 6 & $-0.726\cdot 10^{2}$ & 1.4 & 33/35 & \\
 REA2 & 4 & 4 & 4 &$-0.603\cdot 10^{2}$ & 1.3 & 34/58 & \\
 DIS1 & 8 & 8 & 16 & $-0.117\cdot 10^{2}$ & 137.6 & 30/55 &  \\
 DIS2 & 4 & 3 & 4 & $-0.640\cdot 10^{1}$ & 1.6 & 59/84 & \\
 DIS3 & 8 & 6 & 16 & $-0.168\cdot 10^{2}$ & 642.3 & 66/102 & \\
 MFP & 3 & 4 & 6 & $-0.370\cdot 10^{-1}$ & 1.0 & 20/21 & \\
 TF1 & 4 & 7 & 8 & $-0.847\cdot 10^{-8}$ & 1.7 & 27/31 & 0 \\
 TF2 & 4 & 7 & 6 &  $-0.949\cdot 10^{-7}$ & 1.3 & 19/23 & 0 \\
 TF3 & 4 & 7 & 6 &  $-0.847\cdot 10^{-8}$ & 1.6 & 28/38 & 0 \\
 PSM & 4 & 7 & 6 &  $-0.731\cdot 10^{2}$ & 1.1 & 17/39 & \\
 NN1 & 2 & 3 & 2 &  $-0.131\cdot 10^{0}$ & 1.2 & 32/34 & 0\\
 NN3 & 2 & 4 & 1 &  $0.263\cdot 10^{2}$ & 1.0 & 31/36 & +\\
 NN4 & 4 & 4 & 6 & $-0.187\cdot 10^{2}$ & 1.2 & 33/47 & \\
 NN5 & 2 & 7 & 2 & $0.137\cdot 10^{2}$ & 1.5 & 108/118 & +\\
 NN8 & 3 & 3 & 4 &  $-0.103\cdot 10^{1}$ & 1.0 & 19/29 & \\
 NN9 & 4 & 5 & 6 & $0.312\cdot 10^{1}$ & 1.6 & 64/97 & +\\
 NN10 & 6 & 8 & 9 & $0.409\cdot 10^{4}$ & 18.3 &300/543 & F\\
 NN12 & 4 & 6 & 4 & $0.473\cdot 10^{1}$ & 1.4 & 47/58 & + \\
 NN13 & 4 & 6 & 4 & $0.279\cdot 10^{12}$ & 2.2 & 200/382 &F\\
 NN14 & 4 & 6 & 4 & $0.277\cdot 10^{12}$ & 2.3 & 200/382 &F\\
 NN15 & 3 & 3 & 4 & $-0.226\cdot 10^{0}$ & 1.0 & 15/14 & \\
 NN16 & 7 & 8 & 16 & $-0.623\cdot 10^{3}$ &613.3 &111/191 &  \\
 NN17 & 2 & 3 & 2 & $0.931\cdot 10^{-1}$ & 1.0 & 25/26 & +\\
\hline
  \end{tabular*}         \label{tab:11}
\end{table}
The reader can see that PENLAB can solve all problems apart from AC7,
NN10, NN13 and NN14; these problems are, however, known to be very
ill-conditioned and could not be solved via the Lyapunov matrix
approach either (see \cite{SOF2004}). Notice that the largest problems
with polynomials of degree up to 8 did not cause any major difficulties
to the algorithm.

\paragraph{PENLAB distribution}
The related MATLAB programs are stored in directory \verb|applications/SOF| of the
PENLAB distribution. To solve, for instance, example AC1, run
\begin{verbatim}
>> sof('AC1');
\end{verbatim}
COMPlib program and library must be installed on user's computer.

\section{PENLAB versus PENNON (MATLAB versus C)}\label{sec:comparison}
The obvious concern of any user will be, how fast (or better, how slow)
is the MATLAB implementation and if it can solve any problems of
non-trivial size. The purpose of this section is to give a very rough
comparison of PENLAB and PENNON, i.e., the MATLAB and C implementation
of the same algorithm. The reader should, however, not make any serious
conclusion from the tables below, for the following reasons:
\begin{itemize}
\item Both implementations slightly differ. This can be seen on the
    different numbers of iterations needed to solve single
    examples.
\item The difference in CPU timing very much depends on the type of
    the problem. For instance, some problems require
    multiplications of sparse matrices with dense ones---in this
    case, the C implementation will be much faster. On the other
    hand, for some problems most of the CPU time is spent in the
    dense Cholesky factorization which, in both implementations,
    relies on LAPACK routines and thus the running time may be
    comparable.
\item The problems were solved using an Intel i7 processor with two
    cores. The MATLAB implementation used both cores to perform
    \emph{some} commands, while the C implementation only used one
    core. This is clearly seen, e.g., example lame\_emd10 in
    Table~\ref{tab:12}.
\item For certain problems (such as mater2 in Table~\ref{tab:14}),
    most of the CPU time of PENLAB is spent in the user defined
    routine for gradient evaluation. For linear SDP, this only
    amounts to reading the data matrices, in our implementation
    elements of a two-dimensional cell array, from memory. Clearly,
    a more sophisticated implementation would improve the timing.
\end{itemize}
For all calculations, we have used a notebook running Windows 7 (32
bit) on Intel Core i7 CPU M620@2.67GHz with 4GB memory and MATLAB
7.7.0.

\subsection{Nonlinear programming problems}
We first solved selected examples from the COPS collection \cite{cops}
using AMPL interface. These are medium size examples mostly coming from
finite element discretization of optimization problems with PDE
constraints. Table~\ref{tab:12} presents the results.
\begin{table}[htbp]
\caption{Selected COPS examples. CPU time is given in seconds.
Iteration count gives the number of the global iterations in
Algorithm~\ref{algo:1} and the total number of steps of the Newton
method.}
\begin{tabular*}{\hsize}{@{\extracolsep{\fill}}crrrrrrr} \hline
problem & vars& constr. & constraint & \multicolumn{2}{c}{PENNON} & \multicolumn{2}{c}{PENLAB} \\
        &          &             &       type &  CPU & iter. & CPU & iter. \\\hline
elec200    &  600 &    200 &    $=$ &      40 & 81/224 & 31 & 43/135 \\
chain800   & 3199 &   2400 &    $=$ &       1 & 14/23  &  6 & 24/56\\
pinene400  & 8000 &   7995 &    $=$ &       1 &  7/7   & 11 & 17/17\\
channel800 & 6398 &   6398 &    $=$ &       3 &  3/3   &  1 & 3/3\\
torsion100 & 5000 &  10000 &   $\leq$   &   1 & 17/17  & 17 & 26/26 \\
bearing100 & 5000 &   5000 &   $\leq$   &   1 & 17/17  & 13 & 36/36 \\
lane\_emd10& 4811 &     21 &   $\leq$   & 217 & 30/86  & 64 & 25/49\\
dirichlet10& 4491 &     21 &   $\leq$   & 151 & 33/71  & 73 & 32/68 \\
henon10    & 2701 &     21 &   $\leq$   &  57 & 49/128 & 63 & 76/158 \\
minsurf100 & 5000 &   5000 &   box      &   1 & 20/20  & 97 & 203/203 \\
gasoil400  & 4001 &   3998 &  $=$ \& box &   3 & 34/34  & 13 & 59/71\\
duct15     & 2895 &   8601 & $=$ \& $\leq$&   6 & 19/19  &  9 & 11/11\\
tri\_turtle& 3578 &   3968 &$\leq$ \& box&   3 & 49/49  &  4 & 17/17\\
marine400  & 6415 &   6392 &$\leq$ \& box&   2 & 39/39  & 22 & 35/35 \\
steering800& 3999 &   3200 &$\leq$ \& box&   1 &  9/9   &  7 & 19/40 \\
methanol400& 4802 &   4797 &$\leq$ \& box&   2 & 24/24  & 16 & 47/67 \\
catmix400  & 4398 &   3198 &$\leq$ \& box&   2 & 59/61  & 15 & 44/44 \\
\hline
\end{tabular*}    \label{tab:12}
\end{table}

\subsection{Linear semidefinite programming problems}
We solved selected problems from the SDPLIB collection
(Table~\ref{tab:13}) and Topology Optimization collection
(Table~\ref{tab:14}); see \cite{sdplib,topo}. The data of all problems
were stored in SDPA input files \cite{sdpa}. Instead of PENNON, we have
used its clone PENSDP that directly reads the SDPA files and thus avoid
repeated calls of the call back functions. The difference between
PENNON and PENSDP (in favour of PENSDP) would only be significant in
the mater2 example with many small matrix constraints.
\begin{table}[htbp]
\caption{Selected SDPLIB examples. CPU time is given in seconds.
Iteration count gives the number of the global iterations in
Algorithm~\ref{algo:1} and the total number of steps of the Newton
method.}
\begin{tabular*}{\hsize}{@{\extracolsep{\fill}}crrrrrrr} \hline
problem & vars& constr. & constr. & \multicolumn{2}{c}{PENSDP} & \multicolumn{2}{c}{PENLAB} \\
        &          &             &       size &  CPU & iter. & CPU & iter. \\\hline
control3    &  136 & 2 & 30 &      1 & 19/103 & 20 & 22/315 \\
maxG11  &  800 &  1&  1600 &      18 & 22/41 & 186 & 18/61 \\
qpG11    &  800 &  1&  1600 &      43 & 22/43 & 602 & 18/64 \\
ss30    &  132 &  1&  294 &      20 & 23/112 & 17 & 12/63 \\
theta3    &  1106 &  1&  150 &      11 & 15/52 & 61 & 14/48 \\
\hline
\end{tabular*}    \label{tab:13}
\end{table}

\begin{table}[htbp]
\caption{Selected TOPO examples. CPU time is given in seconds.
Iteration count gives the number of the global iterations in
Algorithm~\ref{algo:1} and the total number of steps of the Newton
method.}
\begin{tabular*}{\hsize}{@{\extracolsep{\fill}}crrrrrrr} \hline
problem & vars& constr. & constr. & \multicolumn{2}{c}{PENSDP} & \multicolumn{2}{c}{PENLAB} \\
        &          &             &       size &  CPU & iter. & CPU & iter. \\\hline
buck2    &  144 & 2 & 97 &      2 & 23/74 & 22 & 18/184 \\
vibra2  &  144 &  2&  97 &      2 & 34/132 & 35 & 20/304 \\
shmup2  &  200 &  2&  441 &     65 & 24/99 & 172 & 26/179 \\
mater2    &  423 & 94 & 11 &      2 & 20/89 & 70 & 12/179 \\
\hline
\end{tabular*}    \label{tab:14}
\end{table}

\appendix
\section{Appendix: Differential calculus for functions of symmetric matrices}
Matrix differential calculus---derivatives of functions depending on
matrices---is a topic covered in several papers; see, e.g.,
\cite{boik2006lecture,dwyer1948symbolic,matrixcookbook,pollock1985tensor}
and the book \cite{magnus1988matrix}. The notation and the very
definition of the derivative differ in these papers. Hence, for
reader's convenience, we will give a basic overview of the calculus for
some typical (in semidefinite optimization) functions of matrices.

For a matrix $X$ (whether symmetric or not), let $x_{ij}$ denote its
$(i,j)$-th element. Let further $E_{ij}$ denote a matrix with all
elements zero except for a unit element in the $i$-the row and $j$-th
column (the dimension of $E_{ij}$ will be always clear from the
context). Our differential formulas are based on Definitions
\ref{def:1} and \ref{def:2}, hence we only need to find the partial
derivative of a function $F(X)$, whether matrix or scalar valued, with
respect to a single element $x_{ij}$ of $X$.

\subsection{Matrix valued functions} Let $F$ be a differentiable
$m\times n$ real matrix function of an $p\times q$ matrix of real
variables $X$. Table~\ref{tab:a1} gives partial derivatives of $F(X)$
with respect to $x_{ij}$, $i=1,\ldots,p$, $j=1,\ldots,q$ for some most
common functions. In this table, $E_{ij}$ is always of the same
dimension as $X$.
\begin{table}[htp]
\caption{}
\begin{tabular}{lcl}
 \hline \\[-1mm]
     $F(X)$  & $\frac{\partial F(X)}{\partial x_{ij}}$ & Conditions\\   [2mm]
 \hline    \\[-1mm]
 $X$&  $E_{ij}$ \\
 $X^T$&  $E_{ji}$ \\
 $AX$&  $AE_{ij}$ &$A\in\RR^{m\times p}$\\
 $XA$&  $E_{ij}A$ &$A\in\RR^{m\times p}$\\
 $XX$& $E_{ij}X + XE_{ij}$ & \\
 $X^TX$& $E_{ji}X + X^TE_{ij}$ & \\
 $XX^T$& $E_{ij}X^T + XE_{ji}$ & \\
 $X^s$& $\displaystyle E_{ij}X^{s-1}+\sum\limits_{k=1}^{s-2} X^kE_{ij}X^{s-k-1} + X^{s-1}E_{ij}$ & $X$ square, $p=1,2,\ldots$\\
 $X^{-1}$ & $-X^{-1} E_{ij} X^{-1}$ & $X$ nonsingular\\
 \hline
\end{tabular}\label{tab:a1}
\end{table}
To compute other derivatives, we may use the following result on the
chain rule.
\begin{theorem}\label{th:a1}
Let $F$ be a differentiable $m\times n$ real matrix function of an
$p\times q$ matrix $Y$ that itself is a differentiable function $G$ of
an $s\times t$ matrix of real variables $X$, that is $F(Y)=F(G(X))$.
Then
\begin{equation}\label{eq:chain}
  \frac{\partial F(G(X))}{\partial x_{ij}} = \sum_{k=1}^p\sum_{\ell=1}^q
  \frac{\partial F(Y)}{\partial y_{k\ell}}\frac{\partial [G(X)]_{k\ell}}{\partial x_{ij}}\,.
\end{equation}
\end{theorem}
In particular, we have
\begin{align}
  \frac{\partial (G(X)H(X))}{\partial x_{ij}} &
  = \frac{\partial G(X)}{\partial x_{ij}}H(X) + G(X)\frac{\partial H(X)}{\partial x_{ij}}\\
  \frac{\partial (G(X))^{-1}}{\partial x_{ij}} &
  = -(G(X))^{-1}\frac{\partial G(X)}{\partial x_{ij}}(G(X))^{-1}\,.
\end{align}
We finish this section with the all important theorem on derivatives of
functions of \emph{symmetric} matrices.
\begin{theorem}\label{th:a2}
Let $F$ be a differentiable $n\times n$ real matrix function of a
symmetric $m\times m$ matrix of real variables $X$. Denote $Z_{ij}$ be
the $(i,j)$-th element of the gradient of $F(X)$ computed by the
general formulas in Table~\ref{tab:a1} and Theorem~\ref{th:a1}. Then
$$
  [\nabla F(X)]_{i,i} = Z_{ii}
$$
and
$$
  [\nabla F(X)]_{i,j} = Z_{ij} + Z_{ji}\qquad \mbox{for $i\not=j$}\,.
$$
\end{theorem}

\paragraph{Example}
Let $X=\begin{pmatrix}x_{11}&x_{12}\\x_{21}&x_{22}\end{pmatrix}$ and
$F(X) = X^2 =
\begin{pmatrix}x_{11}^2+x_{12}x_{21}&x_{11}x_{12}+x_{12}x_{22}\\
x_{11}x_{21}+x_{21}x_{22}&x_{12}x_{21}+x_{22}^2\end{pmatrix}$. Then
$$
  \nabla F(X) = \left[\begin{array}{cc}
  \begin{pmatrix} 2x_{11}& x_{12}\\x_{21} & 0     \end{pmatrix} &
  \begin{pmatrix} x_{21}&x_{11}+x_{22}\\0     &x_{21}\end{pmatrix} \\
  \begin{pmatrix} x_{12}       &  0\\x_{11}+x_{22}&x_{12}\end{pmatrix} &
  \begin{pmatrix} 0     &x_{12}\\x_{21}&2x_{22}\end{pmatrix} \end{array}\right]
$$
(a $2\times 2$ array of $2\times 2$ matrices). If we now assume that
$X$ is symmetric, i.e. $x_{12}=x_{21}$, we get
$$
  \nabla F(X) = \left[\begin{array}{cccc}
  \begin{pmatrix} 2x_{11}& x_{21}\\x_{21} & 0     \end{pmatrix} &
\begin{pmatrix} 2x_{21}       &  x_{11}+x_{22}\\x_{11}+x_{22}&2x_{21}\end{pmatrix} \\
  \begin{pmatrix} 2x_{21}       &  x_{11}+x_{22}\\x_{11}+x_{22}&2x_{21}\end{pmatrix} &
  \begin{pmatrix} 0     &x_{21}\\x_{21}&2x_{22}\end{pmatrix} \end{array}\right]\,.
$$
We can see that we could obtain the gradient for the symmetric matrix
using the general formula in Table~\ref{tab:a1} together with
Theorem~\ref{th:a2}.

Notice that if we simply replaced each $x_{12}$ in $\nabla F(X)$ by
$x_{21}$ (assuming symmetry of $X$), we would get an \emph{incorrect}
result
$$
  \nabla F(X) = \left[\begin{array}{cc}
  \begin{pmatrix} 2x_{11}& x_{21}\\x_{21} & 0     \end{pmatrix} &
  \begin{pmatrix} x_{21}&x_{11}+x_{22}\\0     &x_{21}\end{pmatrix} \\
  \begin{pmatrix} x_{21}       &  0\\x_{11}+x_{22}&x_{21}\end{pmatrix} &
  \begin{pmatrix} 0     &x_{21}\\x_{21}&2x_{22}\end{pmatrix} \end{array}\right]\,.
$$

\subsection{Scalar valued functions}

Table~\ref{tab:a5} shows derivatives of some most common scalar valued
functions of an
 $m\times n$ matrix~$X$.
\begin{table}[htp]
\caption{}
\begin{tabular}{lccl}
 \hline              \\[-1mm]
     $F(X)$ & equivalently  & $\frac{\partial F(X)}{\partial x_{ij}}$ & Conditions\\  [2mm]
 \hline      \\[-1mm]
 $\Tr X$ & $\langle I,X\rangle$ &  $\delta_{ij}$ \\
 $\Tr AX^T$ & $\langle A,X\rangle$ &  $A_{i,j}$ &$A\in\RR^{m\times n}$ \\
  $a^TXa$ & $\langle aa^T,X\rangle$ &  $a_ia_j$ &$a\in\RR^n$, $m=n$\\
 $\Tr X^2$ & $\langle X,X\rangle$ &  $2X_{j,i}$ &$m=n$\\
 \hline
\end{tabular}\label{tab:a5}
\end{table}

Let $\Phi$ and $\Psi$ be functions of a square matrix variable $X$. The
following derivatives of composite functions allow us to treat many
practical problems (Table~\ref{tab:a6}).
\begin{table}[htp]
\caption{}
\begin{tabular}{lccl}
 \hline \\[-1mm]
     $F(X)$ & equivalently  & $\frac{\partial F(X)}{\partial x_{ij}}$ & Conditions\\   [2mm]
 \hline \\[-1mm]
 $\Tr A\Phi(X)$ & $\langle A,\Phi(X)\rangle$ &  $\langle A, \frac{\partial \Phi(X)}{\partial x_{ij}}\rangle$ \\
 $\Tr \Phi(X)^2$ & $\langle \Phi(X),\Phi(X)\rangle$ &  $2\langle \Phi(X), \frac{\partial \Phi(X)}{\partial x_{ij}}\rangle$ \\
 $\Tr (\Phi(X)\Psi(X))$ & $\langle \Phi(X),\Psi(X)\rangle$ &  $\langle \Phi(X), \frac{\partial \Psi(X)}{\partial x_{ij}}\rangle
 +\langle \frac{\partial \Phi(X)}{\partial x_{ij}},\Psi(X)\rangle$ \\
 \hline
\end{tabular}\label{tab:a6}
\end{table}
We can use it, for instance, to get the following two results for a
$n\times n$ matrix $X$ and $a\in\RR^n$:
\begin{align*}
  \frac{\partial}{\partial x_{ij}} (a^TX^{-1}a) &= \frac{\partial}{\partial x_{ij}}\langle aa^T,X^{-1}\rangle\\
  &=-\langle aa^T,X^{-1} E_{ij}X^{-1}\rangle \\
  & = -a^TX^{-1}E_{ij}X^{-1}a\,,
\end{align*}
in particular,
$$
  \frac{\partial}{\partial x_{ij}} \Tr X^{-1} = \frac{\partial}{\partial x_{ij}}\langle I,X^{-1}\rangle =
   -\langle I,X^{-1} E_{ij}X^{-1}\rangle = - \Tr (X^{-1}E_{ij}X^{-1})\,.
$$
Recall that for a \emph{symmetric} $n\times n$ matrix $X$, the above
two formulas would change to
$$
  \frac{\partial}{\partial x_{ij}} (a^TX^{-1}a)
  = -a^T(Z_{ij}+Z_{ij}^T-{\rm diag}Z_{ij})a
$$
and
$$
  \frac{\partial}{\partial x_{ij}} \Tr X^{-1}
  = - \Tr (Z_{ij}+Z_{ij}^T-{\rm diag}Z_{ij})
$$
with
$$
  Z_{ij} = X^{-1}E_{ij}X^{-1}\,.
$$

\subsection{Second-order derivatives}
To compute the second-order derivatives of functions of matrices, we
can simply apply the formulas derived in the previous sections to the
elements of the gradients. Thus we get, for instance,
$$
  \frac{\partial^2 X^2}{\partial x_{ij}\partial x_{k\ell}} =
  \frac{\partial}{\partial x_{k\ell}} (E_{ij}X + XE_{ij}) =
  E_{ij}E_{k\ell} + E_{k\ell}E_{ij}
$$
or
$$
  \frac{\partial^2 X^{-1}}{\partial x_{ij}\partial x_{k\ell}}
  =\frac{\partial}{\partial x_{k\ell}} (-X^{-1}E_{ij}X^{-1})
  =X^{-1}E_{ij}X^{-1}E_{kl}X^{-1} + X^{-1}E_{kl}X^{-1}E_{ij}X^{-1}
$$
for the matrix valued functions, and
$$
   \frac{\partial^2}{\partial x_{ij}\partial x_{k\ell}}\langle \Phi(X),\Phi(X)\rangle
   = 2\left\langle \frac{\partial \Phi(X)}{\partial x_{kl}},\frac{\partial \Phi(X)}{\partial x_{ij}}\right\rangle
    +2\left\langle \Phi(X),\frac{\partial^2 \Phi(X)}{\partial x_{ij}\partial x_{k\ell}}\right\rangle
$$
for scalar valued matrix functions. Other formulas easily follow.

\bibliographystyle{plain}
\bibliography{pennon,sensor}

\end{document}